\documentclass[10pt]{amsart} 
\RequirePackage{graphicx}

\usepackage[top=3cm,bottom=2cm,left=3cm,right=3cm,marginparwidth=1.75cm]
{geometry}
\usepackage{mathrsfs}   
\usepackage[english]{babel}
\usepackage[T1]{fontenc}
\usepackage{xcolor, import}

\usepackage{dsfont}

\usepackage{fancyvrb}
\usepackage[colorlinks=true, citecolor=red, linkcolor=blue]{hyperref}
\usepackage{xparse}
\usepackage[titletoc,title]{appendix}
\usepackage{amsfonts,mathtools}
\usepackage{amsmath,amssymb,amsthm}
\usepackage{bbm,xargs}
\usepackage[square,compress,comma, numbers,sort]{natbib}
\usepackage{mathbbol}
\usepackage[shortlabels]{enumitem}
\usepackage[nameinlink,noabbrev]{cleveref}

\usepackage{aliascnt}
 \theoremstyle{plain}
\newtheorem{theorem}{Theorem}[section]

\newaliascnt{lemma}{theorem}
\newtheorem{lemma}[lemma]{Lemma}
\aliascntresetthe{lemma}

\newaliascnt{proposition}{theorem}

\aliascntresetthe{proposition}

\newaliascnt{corollary}{theorem}
\newtheorem{corollary}[corollary]{Corollary}
\aliascntresetthe{corollary}

\theoremstyle{definition}
\newaliascnt{definition}{theorem}
\newtheorem{definition}[definition]{Definition}
\aliascntresetthe{definition}

\newaliascnt{example}{theorem}
\newtheorem{example}[example]{Example}
\aliascntresetthe{example}

\theoremstyle{remark}
\newaliascnt{remark}{theorem}
\newtheorem{remark}[remark]{Remark}
\aliascntresetthe{remark}

\crefname{theorem}{Theorem}{Theorems}
\Crefname{theorem}{Theorem}{Theorems}
\crefname{lemma}{Lemma}{Lemmas}
\Crefname{lemma}{Lemma}{Lemmas}
\crefname{proposition}{Proposition}{Propositions}
\Crefname{proposition}{Proposition}{Propositions}
\crefname{corollary}{Corollary}{Corollaries}
\Crefname{corollary}{Corollary}{Corollaries}
\crefname{definition}{Definition}{Definitions}
\Crefname{definition}{Definition}{Definitions}
\crefname{example}{Example}{Examples}
\Crefname{example}{Example}{Examples}
\crefname{remark}{Remark}{Remarks}
\Crefname{remark}{Remark}{Remarks}

\definecolor{c20}{rgb}{0.,0.7,0.}
\definecolor{c30}{rgb}{0.,0.,1.}
\definecolor{c40}{rgb}{1,0.1,0.7}
\definecolor{c50}{rgb}{1,0,0}
\definecolor{c60}{rgb}{1,0.9,0.1}
\definecolor{c70}{rgb}{0.50,1.00,0.00}

\newcommand{\R}{\mathbb{R}}
\newcommand{\inr}{\in \R}
\newcommand{\QED}{\hfill $\Box$}

\DeclarePairedDelimiter{\parens}()
\def\Pspace{ \parens{\Omega,\mathscr{F},\mathbb{P}}}

\DeclarePairedDelimiterXPP\pk[1]{\mathbb{P}}\{ \}{}{ #1}
\DeclarePairedDelimiterXPP\pkh[1]{{\mathbb{P}}}\{ \}{}{ #1}
\DeclarePairedDelimiterXPP\pkT[1]{{\mathbb{P}}}\{ \}{}{ #1}
\DeclarePairedDelimiterXPP\pTI[1]{{\mathbb{P}}}\{ \}{}{ #1}
\DeclarePairedDelimiterXPP\E[1]{\mathbb{E}}\{ \}{}{	#1}
\DeclarePairedDelimiterXPP\EHh[1]{{\mathbb{E}}}\{ \}{}{	#1}
\DeclarePairedDelimiterXPP\ETI[1]{{\mathbb{E}}}\{ \}{}{	#1}

\def\kk{\mathcal{C}  [\vk Z ]}

\def\IF{\infty}

\newcommand{\kb}[1]{\boldsymbol{#1}}
\newcommand{\vk}[1]{\kb{#1}}

\newcommand{\COM}[1]{}

\newcommand\Dset{(\R^d)^\TT}

\def\Hh{ \mathfrak{H}}

\def\clasP{\mathfrak{W}}

\DeclarePairedDelimiterXPP\ind[1]{\mathbb{I}}( ){}{	#1}

\def\toas{\overset{a.s.}\rightarrow}
\def\toprob{\overset{p}\rightarrow}

\newcommand{\norm}[1]{\lVert  #1 \rVert }

\newcommand{\Nset}{\mathbb{N}}

\def\boldsymbol#1{#1}

\def\AA{\mathscr{A}}

\newcommand{\prooftheo}[1]{ \textsc{\bf Proof of Theorem} \ref{#1}: }

\newcommand{\prooflem}[1]{\textsc{\bf Proof of Lemma} \ref{#1}:}
\newcommand{\proofkorr}[1]{\textsc{\bf Proof of Corollary} \ref{#1}:}

\newcommand{\BEL}{\begin{lemma}}
\newcommand{\EEL}{\end{lemma}}

\newcommand{\BRM}{\begin{remark}}
\newcommand{\ERM}{\end{remark}}

\def\TT{\mathcal{T} }
\def\intT{\int_{\TT}}
\def\TT{\mathcal{T} }
\def\ZT{\widetilde{ \vk Z}}

\def\diad{{\mathbb{T}_0}}

\def\bqny#1{ \begin{eqnarray*} #1 \end{eqnarray*}}
\def\bqn#1{ \begin{eqnarray} #1 \end{eqnarray}}
  
\newcommand{\BT}{\begin{theorem}}
\newcommand{\ET}{\end{theorem}}
\newcommand{\BK}{\begin{corollary}}
\newcommand{\EK}{\end{corollary}}
\newcommand{\BD}{\begin{definition}}
\newcommand{\ED}{\end{definition}}

\newcommand{\BEX}{\begin{example}}
\newcommand{\EEX}{\end{example}}
\def\bex#1{\begin{example}{\it #1} \end{example}}

\def\TT{\mathcal{T}}

\def\diad{{\mathbb{T}_0}}

\def\SZT{ \mathcal{S}( \ZT)} 
\def\supA{\sup_{t\in \diad}}
\def\supA{\sup_{t\in \TT}}

\begin{document}

\title[Shift-generated  classes ]
{Shift-generated   classes of jointly measurable random fields}

	\author{Enkelejd  Hashorva}
	\address{Enkelejd Hashorva, HEC Lausanne,  
		University of Lausanne,\\
		Chamberonne, 1015 Lausanne, Switzerland}
	\email{Enkelejd.Hashorva@unil.ch}

	\bigskip
	
	\date{\today}
	\maketitle

	\begin{quote} 
		{\bf Abstract:} We study shift-generated classes of jointly measurable and separable
\(\mathbb R^d\)-valued random fields (RFs) indexed by \(\mathbb R^l\), defined
through identities for \(\alpha\)-homogeneous functionals. In contrast to
earlier work, no stochastic-continuity assumption and no local boundedness
condition are imposed. We show that every non-empty shift-generated class
contains an \(L^\alpha\)-continuous element. This regularization result
allows us to establish the strict positivity of the integral functional  for all elements of the class and for the associated local
RFs. We further extend the defining functional identity to a
larger class of functionals, including integral functionals, and use this to
construct canonical elements of a given class via randomised shifts. We also
relate shift-generated classes to spectral tail and tail RFs and
show that every spectral tail RF has an \(L^\alpha\)-continuous
representative with the same finite-dimensional distributions. As an
application, we identify the \(-\alpha\)-homogeneous tail measure associated
with a shift-generated class and show that it depends only on the class and
admits an \(L^\alpha\)-continuous representor.
 	\end{quote}

    

\section{Introduction}
Let $\clasP$ denote the class of jointly measurable and separable
$\mathbb R^d$-valued RFs $U(t),t\in\TT$, where
$\TT=\mathbb R^l$ and $d,l\in\mathbb N$, defined on a complete non-atomic
probability space $\Pspace$.
 
Given a fixed positive constant $\alpha>0$, 
let  the RF $Z\in \clasP$  satisfy
\bqn{\label{eq1}
	\pk*{ \supA  \norm{Z(t)}>0 }=1, \quad \E{\norm{Z(0)}^\alpha}=C\in (0,\IF),  
}	 
where $\norm{\cdot}$ is some norm on $\R^d$.
  
Equip $\Dset$ with the product $\sigma$-field $\AA$ and   
let $\mathcal H$ be the class of all $\AA/\mathcal B([0,\infty])$-measurable
maps $F:\Dset\to[0,\infty]$.     Write $\mathcal{H}_\beta$ for the subset of $\mathcal{H}$ consisting of $\beta$-homogeneous maps $F$,  i.e.,  
$F(c f) = c^\beta F( f)$ for all $f \in \Dset,c>  0$.\\

In this contribution we are interested in the class of RFs $\ZT \in \clasP$, which is {\it shift-generated} from the RF $Z$ through the  {\it functional identity} \eqref{boll} below. For the special case that 
\bqn{\label{specC} 
\E*{\sup_{t\in [-a,a]^l} \norm{Z(t)}^\alpha}<\infty, \qquad \forall 
a>0}
these classes were 
considered in \cite{hashorva2021shiftinvariant}.  In the following $B^h f(t)=f(t-h), h\in \TT, f\in \Dset$ denotes the shift operator.

\BD
Let \(Z\in\clasP\) satisfy \eqref{eq1} and let  the class of RFs 
$\kk$ consist of all $\ZT\in\clasP$ that satisfy \eqref{eq1} with the same constant $C$,  
such that for every $h\in\TT$ 
\bqn{\label{boll}
		\E*{ F(Z)}=\ETI{ F(B^h\ZT)},\qquad \forall F\in\mathcal H_\alpha.
}   
\ED  
If $\kk$ is non-empty, and thus there exists $\ZT \in \kk$, by the definition above we have $B^h \ZT \in \kk$ for all $h\in \TT$, which justifies the term {\it shift-generated}.

For a given stationary  $Z\in \clasP$  such that   \eqref{eq1} holds, the corresponding class $\kk$ is clearly non-empty since $Z\in \kk$.  
Not every RF $Z\in \clasP$ generates a non-empty class $\kk$.  We give next an example, known from \cite{kab2009}, which shows non-stationary $Z$'s that define non-empty  $\kk$'s.

\bex{(Brown--Resnick shift-generated RFs)
For \(\alpha=1\), consider the log-normal RF
\[
        Z_i(t)=\exp(W_i(t)-\operatorname{Var}(W_i(t))/2),
        \qquad i=1,\ldots,d,\quad t\in\TT,
\]
where \(W(t)=(W_1(t),\ldots,W_d(t))\) is a centered
\(\mathbb R^d\)-valued Gaussian RF. Define the  
pseudo-variogram $\gamma=(\gamma_{ij}),1\le i,j\le d$ of $W$ by
\[
        \gamma_{ij}(s,t)
        =
        \operatorname{Var}(W_i(s)-W_j(t)),
        \qquad 1\le i,j\le d .
\]
Then the law of the associated Brown--Resnick max-stable RF depends only on
\(\gamma=(\gamma_{ij})\). In particular, if
\[
        \gamma_{ij}(s+h,t+h)=\gamma_{ij}(s,t),
        \qquad s,t,h\in\TT,\quad 1\le i,j\le d,
\]
equivalently if \(\gamma_{ij}(s,t)\) depends only on \(t-s\), then \(X_Z\)
is stationary and \(Z\) generates a non-empty shift-generated class
\(\mathcal C[Z]\), see \cite[Lem 4.2]{KumeE}. Conversely, stationarity of the 
Brown--Resnick RF
forces the same translation-invariance of the pairwise pseudo-variograms
\(\gamma_{ij}\).
}
Given two independent $Z_1,Z_2\in \clasP$ satisfying \eqref{eq1} with the same constant $C$, the class  generated by $Z_1*Z_2$ (component-wise/Hadamard multiplication) is non-empty if  both $\mathcal{C}[Z_1],\mathcal{C}[Z_2]$ are non-empty. An interesting model introduced in \cite{hashorva2021shiftinvariant} is by taking $Z_2(t)=\xi$ a $d$-dimensional random vector independent of $Z_1$.\\ 
 Another possibility to generate a new class is by defining $Z_I=Z_1$ if $I=1$ and $Z_I=Z_2$ if $I$ equals 2, with $I$ a random variable (rv) taking values in $\{1, 2\}$  being independent of $Z_1$ and $Z_2$.  \\

 In the particular case that  $Z $ has non-negative components, the class $\kk$ is closely related to the class of max-stable RFs $X_Z$, which are defined via the de Haan representation (\cite{deHaan,dom2016}) 
\bqn{\label{eqDehaan}
	{X}_Z (t)=  \max_{i\ge 1}  \frac{ Z^{(i)}(t)}{   (\sum_{k=1}^i \mathcal{V}_k)^{1/\alpha}} , \quad t\in \TT,
}  
where   $Z^{(i)}$'s are 
  independent copies of $Z$ 
 being further independent of and identically distributed (iid) unit exponential rvs $
\mathcal{V}_k, k\ge 1$. 
Here the   maximum is applied component-wise and $Z$ is called a representor of $X_Z$. 

 The importance of a given  
   non-empty $\kk$ relates to the fact that   $X_{\ZT}$ for some $\ZT \in \kk$   
 is stationary and has the same law as  $ X_{Z}$ and conversely, if $X_Z$ is stationary, then $Z$ generates a non-empty $ \kk$, see e.g.,  
 see \cite{MolchanovBE,Htilt, MartinE, kulik:soulier:2020,  hashorva2021shiftinvariant}.  
 
 For general $Z\in \kk$ one can further define a corresponding stationary stable process, see e.g., \cite{PH2020, kulik:soulier:2020}.

The class $\kk$ was   introduced in \cite{hashorva2021shiftinvariant} and further investigated in \cite{hashorva2025cluster}. Therein it was assumed that $Z$ is   stochastically continuous  and $X_Z$ has 
locally bounded  sample paths; both assumptions are dropped in this contribution. 

Throughout, unless stated otherwise, $\gamma$ denotes a strictly positive
continuous probability density (pdf) on $\TT$, and $N$ denotes an independent
$\TT$-valued random vector with pdf  $\gamma$. We shall also use the
unweighted functional corresponding formally to $\gamma\equiv1$; in this case
$\gamma$ is not a pdf and no random shift $N$ with pdf 
$\gamma$ is involved.

 If the max-stable stationary RF   $X_Z$ 
 is locally bounded,   its representor $Z \in \clasP$ is stochastically
  continuous   and further  almost surely 
\bqn{\label{ataka}
\mathcal{S}(Z)=\intT \norm{Z(t)}^\alpha \lambda(dt)<\IF 
} 
 with   $\lambda(dt)$   the Lebesgue measure on $\TT$, 
 then in view of \cite{hashorva2025cluster} the  max-stable
  RF $ X_{Z_N}$ with representor $Z_N $  given by 
 \bqn{
	Z_{ N}(t)=  \frac{\norm{Z(0)}}{ (\gamma (N)\mathcal{S}(Z))^{1/\alpha }} 
	B^N Z(t) , \quad t\in \TT
	\label{boll5}
}
has the same law as the max-stable RF $X_Z$.  Here the integral functional $\mathcal{S}_\gamma(\cdot)$ is defined by  
$$
\mathcal{S}_\gamma:f \mapsto \intT \norm{f(t)}^\alpha \gamma(t)\lambda(dt),\quad \text{if } f\in \Dset  \text{ is Lebesgue measurable}
$$	
and $\mathcal{S}_\gamma(f)=\IF$ otherwise. 
 We   omit the subscript $\gamma$ writing   simply $\mathcal{S}$, 
 if $\gamma(t)=1,t\in \TT$. \\
 
When condition \eqref{ataka} is not satisfied,  
we still can define 
 \bqny{
	 Z_N'(t)= 
	\norm{Z(0)}\frac{B^N Z(t)}{(\mathcal{S}_\gamma(B^N Z))^{1/\alpha}} ,
	\quad t\in \TT.
}

 By our assumption the RF $Z$ is jointly measurable, in view of   \eqref{eq1} and  \eqref{boll},   
 the Tonelli Theorem yields that $\mathcal{S}_\gamma(Z)$ is a  
 well-defined non-negative rv.  However, for  $Z_N'$ to be well-defined and finite-valued   (and similarly for $Z_N$) we need  to prove  that 
 \bqn{ 
	\label{pandm}
 \pk{\mathcal{S}_\gamma(Z)>0}=1 ,
 }  which for 
 stochastically continuous $Z$ has been shown in 
 \cite{hashorva2021shiftinvariant}.

In order to remove the dependence on $\norm{Z(0)}$,  we shall construct another RF 
$Z_N''$ in terms of the local RF $\widetilde{\Theta}$ of $\kk$. 
Specifically, as in \cite{hashorva2021shiftinvariant}, 
  $\widetilde{ \Theta}$  stands 
for the RF  $\ZT/ \norm{\ZT(0)} $ under the probability measure 
\bqn{ \label{minist}
	\widehat{\mathbb{P}}\{A\}
&=& 
	\ETI{   \norm{  \ZT(0)}^\alpha   \ind{ A} } /\ETI{\norm{  \ZT(0)}^\alpha} , 
	\quad \forall A \in {\mathscr{F}},
}
with $\ind{A}$   the indicator function of some set $A$. We write $\Theta$ if we use the tilting with respect to $Z/\norm{Z(0)}$.
 
In this contribution   
\begin{enumerate}[(i)]
	\item only joint  measurability of   $Z$ will be assumed (dropping stochastic continuity); 
	\item the condition  \eqref{specC} will be dropped. 
\end{enumerate}

In  this  general setting our goals are:
 \begin{enumerate}[C1)]
	
	\item To extend  \eqref{boll} to a large 
    class of maps including $\mathcal{S}_\gamma$, which is crucial to show that 
    $Z_N$ and $Z_N'$ satisfy \eqref{boll};
    \item To show  that a shift-generated $\kk$ has   an $L^\alpha$-continuous 
    element. This will allow us to prove that \eqref{pandm} holds for all $\ZT \in \kk$;
  	\item To identify canonical generators of a given class $\kk$;
  	\item   	 To identify the \(-\alpha\)-homogeneous tail measure associated
    with a shift-generated class and to show that it admits an
    \(L^\alpha\)-continuous representor;
  	\item To  construct  non-empty  $\kk$'s via {\it cluster} RFs. 
   
\end{enumerate}

As for stochastically continuous $Z$,  a 
simple approach to deal with the case  
 $$\pk{ \norm{Z(0)}=0}>0$$ is to introduce the local RFs  $\widetilde{ \Theta}$  as in \eqref{minist}.  
Recall that  $\ZT$ is defined on the complete non-atomic probability 
space $({\Omega}, {\mathcal{F}}, {\mathbb{P}})$. Below we work on the $\widehat {\mathbb{P}}$-completion of $(\Omega,\mathscr F)$. \\
The functional equation \eqref{boll} 
can be written equivalently in terms of 
 $\Theta$, see \eqref{eqDo20} below. 
In turn, such functional identities  can also be utilised 
to define spectral tail RFs, already discussed in 
\cite{Hrovje,klem,kulik:soulier:2020, hashorva2021shiftinvariant} for stochastically continuous 
$\widetilde{\Theta}$. For the latter case, it is known that 
spectral tail RFs define uniquely a $\kk$ and moreover, they appear as limiting RFs 
in the study of heavy-tailed stationary RFs, 
see e.g., \cite{BojanS,BojanPhilippe,kulik:soulier:2020,Resnickart}. 

In this contribution we extend the previously mentioned findings to the larger class of jointly measurable RFs. A further implication concerns homogeneous tail measures, initially appearing in \cite{wao, MR3561100}. To any
\(Z\in\kk\) one may associate the \(-\alpha\)-homogeneous measure obtained by
integrating over radial scalings of \(Z\). The defining identity of
shift-generated classes implies that this measure depends only on the class
\(\kk\), and not on the particular representor. Our regularization result
then shows that this tail measure always admits an \(L^\alpha\)-continuous
representor, showing that  \(L^\alpha\)-regularity is  an intrinsic property of tail measures.\\

The rest of the paper is organised as follows:\\
In Section 2, we present our main results,
 including 
 the validity  of the functional identity \eqref{boll} 
for a large class of functionals and in particular for the integral functional, the existence of $L^\alpha$-continuous elements 
 in shift-generated $\alpha$-homogeneous classes, and the construction of elements of a given $\kk$.  

 Section 3 discusses further implications and extensions, such as the choice
of the norm, the connection to tail and spectral tail RFs, homogeneous tail
measures and their \(L^\alpha\)-regular representatives, the case of locally
bounded sample paths, and the role of cluster RFs for constructing non-empty
shift-generated classes \(\kk\).
 
Section 4 contains the proofs of the main results  followed by an Appendix.

\def\clmapB{\Hh_\beta}
\def\clmapA{\Hh_\alpha}
\def\clmap{\Hh_0}

\section{Main Results} 
In this section we shall consider a non-empty shift-generated   class 
$\kk$ and suppose for simplicity that $\E{\norm{Z(0)}^\alpha}=1$ and hence by  \eqref{boll}  
\bqn{
    \E{\norm{\ZT(t)}^\alpha}=1, \quad \forall t \in \TT, \forall \ZT \in \kk.
 \label{cornidoesta}
}
Below $R$ is a positive rv independent of any other element with survival function
$\widehat{\mathbb{P}}\{R>s\}= s^{-\alpha}$, $s\ge 1$.
Recall the definition of local RFs $\widetilde{\Theta}$ 
in the Introduction   and
write below for simplicity $\mathbb{P},\mathbb{E}$ instead of
$\widehat{\mathbb{P}},\widehat{\mathbb{E}}$. Also note that we write $\Theta$ if we use the tilting with respect to $Z/\norm{Z(0)}$. 
Set
$$
\widetilde{Y}(t)=R \widetilde{\Theta}(t), \quad t \in \TT, \qquad
\mathcal{S}_{\mathbb{I}}(\widetilde{Y})= \intT \ind{\norm{\widetilde{Y}(t)}>1}\,\lambda(dt)
$$
 
 with $\ind{\cdot}$ the indicator function (recall that $\lambda(dt)$ is the Lebesgue measure on $\TT$).  Define similarly $Y$ from $\Theta$ by setting $Y(t)=R \Theta(t), t \in \TT$.

In the definition of $\kk$, the functional identity \eqref{boll} is required
to hold for all \(F\in\mathcal H_\alpha\). We show next that it extends to the
larger class \(\Hh_\alpha\), which contains the integral functional
\(\mathcal S_\gamma\).

\BD
For $\beta\ge0$, let $\Hh_\beta$ be the class of all maps $F:\Dset\to[0,\IF]$ such that
for every $U\in\clasP$ we have that $F(U)$ is a well-defined rv and there exist, for each $n\in\Nset$, points
$t_1^{(n)},\ldots,t_n^{(n)}\in\TT$ and a Borel measurable $\beta$-homogeneous function
$F_n:(\R^d)^n\to[0,\IF]$ such that
\[
F_n\bigl(U(t_1^{(n)}),\ldots,U(t_n^{(n)})\bigr)\toprob F(U), \qquad n\to\IF.
\]
\ED
Note that $\Hh_\beta$ can contain functionals that are not $\AA/\mathcal B([0,\infty])$-measurable,  for instance the integral functional, see \Cref{lemGen2} in Appendix.  
    We extend the functional identity \eqref{boll} to such general $F$'s and show that several interesting functionals including the integral one are contained in $\clmapB$.

\BT 
\label{verbpars}
Given a non-empty shift-generated class $\kk$, we have that  \eqref{boll}  is equivalent to 
 \bqn{ \label{boll2} 
	\E{ F(  Z)}=	\ETI{F(B^h \ZT )}, 
	\quad   \forall h\in \TT, \forall F\in \Hh_\alpha, \forall\ZT \in \kk,
}
which is also equivalent to   
\bqn{ \label{boll22} 
	\E{ \norm{Z(h)}^\alpha G(Z)}=	\ETI{ \norm{\ZT(0)}^\alpha G(B^h \ZT )}, 
	\quad   \forall h\in \TT,\forall G\in \Hh_0, \forall\ZT \in \kk.
}
In particular,
for a given continuous pdf $\gamma(t)>0, t\in \TT$ and $\Gamma_\beta \in \Hh_\beta, \beta \ge 0$, 
 we have 
\begin{enumerate} [(i)]
    \item $F(f)=\Gamma_\beta(f)\ind{\mathcal{S}_\gamma(f) =a}
    \in \Hh_\beta, a\in \{0,\IF\}$;
    \item  $F(f)=\Gamma_\beta (f) \mathcal{S}_\gamma(f)\in 
    \Hh_{\alpha+ \beta}$;
    \item  
    $F(f)=\Gamma_\alpha(f) /\mathcal{S}_\gamma(f)\in 
    \Hh_{0}$.
\end{enumerate}     
\ET

\BRM
\begin{enumerate}[(i)] 
	\item  If $\tilde{Z} \in \clasP$ satisfies \eqref{eq1} and 
	further for all bounded $G\in\Hh_0$ we have 
	$$\E{ \norm{\tilde{Z}(h)}^\alpha G(\tilde{Z}) }
	=\E{G(B^h\Theta)},h\in \TT,$$ then $\tilde{Z} \in \kk$.
	\item A direct implication of \Cref{verbpars}
	 and \cite[Lem 9.7]{hashorva2021shiftinvariant} is that 
	 for any shift-invariant map $F$, i.e., 
	 $F(B^h f)=F(f),f\in \Dset,h\in \TT$ such that $F \in \Hh_0$ we have  
$$ \E{ F(\ZT)}=0 \iff  \E{ F(\widetilde{\Theta})}=0,$$
which implies that $\pk{\mathcal{S}(Z) =a}=0$ is equivalent 
to $\pk{\mathcal{S}(\Theta) =a}=0$ for $a\in \{0,\IF\}$.
\item In view of \cite[Eq. (5.1)]{hashorva2021shiftinvariant},  
if $\kk$ consists of stochastically continuous elements, then  
\eqref{boll} {implies}    for the local RF $\Theta$ defined from $Z$ 
we have
\bqn{\label{eqDo20}
	\E*{
		\norm{\Theta(h)}^\alpha   \Gamma(\Theta)} &=& \ETI{\ind{\norm{\Theta(-h)}\not=0}\Gamma( B^{h}{\Theta}) }, \quad 
\forall \Gamma\in  \Hh_0,\forall h\in \TT.
}
 \Cref{verbpars} shows that the stated identity holds also for the setting of this paper where $\kk$  has jointly measurable elements. 
\end{enumerate}  
\ERM

 We show next  the existence 
 of an $L^\alpha$-continuous  element in 
 $\kk$ and the strict positivity of several rvs of interest. 
\BT 
\label{thmAug23}
There exists an $L^\alpha$-continuous $Z^*\in \kk$ 
and  
\bqn{ \label{pandm2} 
\pk{\mathcal{S}(\ZT )>0}=1, \quad \forall \ZT \in \kk.
}
Moreover, there exists an   $L^\alpha$-continuous  local RF $\Theta^* \in \clasP$ and  for all local RFs $\widetilde{\Theta} \in \clasP$ associated to $\kk$ we have
\bqn{ 
 	\pk{ \mathcal{S}(\widetilde{\Theta})>0}=
    \pk{  \mathcal{S}_{\mathbb{I}}(\widetilde{Y})>0} =1.
\label{e11}
} 
\ET

\BRM \label{sionR}
\begin{enumerate}[(i)] 
\item \label{sionR:1}
All local RFs \(\widetilde\Theta\) of \(\kk\) have the same finite-dimensional
distributions. Hence, by \Cref{lemGen2}, the laws of
\(\mathcal S(\widetilde\Theta)\) and
\(\mathcal S_{\mathbb I}(\widetilde Y)\) do not depend on the choice of
the local RF \(\widetilde\Theta\), respectively of \(\widetilde Y=R\widetilde\Theta\).

\item The claims in \eqref{pandm2} and \eqref{e11} are shown in
\cite[Thm 2.9]{hashorva2025cluster} for stochastically continuous \(Z\)
under the additional assumption \eqref{specC}.
\end{enumerate}
\ERM

The above findings  can now be utilised to construct elements   of $\kk$,  stated in the next result.
 As in the Introduction,  let $N$ be a  
  $\TT$-valued rv being independent of any other random 
  element and with  positive pdf $\gamma(t),t\in \TT$.
 It is defined on the complete non-atomic probability 
 space $(\Omega, \mathcal{F}, \mathbb{P})$ or on 
 $(\Omega, \mathcal{F}, \widehat{\mathbb{P}})$, depending 
 on the context. For simplicity, we shall include in $\kk$ also jointly measurable 
 RFs $\ZT$ defined in the latter probability space which satisfy further 
 \eqref{eq1} and \eqref{boll}.  We consider next only the local RF $\Theta$, but the constructions hold also for $\widetilde {\Theta}$.

\begin{corollary}
	  \label{Theo2}
 
Let $\kk$ with local RF $\Theta$ be given. If $\ZT \in \kk$  and  
$\diad\subset \TT$ is  a countable separant of $\ZT$ and $\Theta$, then 
the following RFs 
    \bqn{
	Z_N(t)= \frac{\norm{\ZT(0)} }{ ( \mathcal{S}_\gamma
     ( B^N \ZT ))^{1/\alpha}}B^N \ZT(t), \qquad  Z_N'(t)= \frac{1}{ (
		\mathcal{S}_\gamma( B^N {\Theta} ))^{1/\alpha}}B^N 
        {\Theta}(t),
	  \quad t \in \TT,
  	\label{boll3}
}
\bqn{
         Z_N''(t)= \frac{B^N 
        {\Theta}(t)}{ 
		( \intT \gamma(s)\ind{\norm{B^N\Theta(s)}>0 } \lambda(ds) )
        ^{1/\alpha}}, \quad t \in \TT
	\label{boll3b}
}
satisfy \eqref{boll}. Moreover, the claims hold 
 also if 
$\lambda(dt)$ 
is  the counting measure on $\diad$ and 
$N$ is  a $\diad$-valued rv with positive pdf $\gamma(t)>0, t\in \diad$.  

 \end{corollary}

 \BRM \label{rLL}
  \begin{enumerate}[(i)] 
        \item In the above construction both $ Z_N'$ and $Z_N''$ belong to       $\kk$, which is also the case for  $Z_N$ when   $\norm{Z(0)}>0$ almost surely. Note that $Z_N$ appears already in \cite{KumeE}, however for the general setting of this paper, the proof therein can be made rigorous only after using our new result \eqref{boll22}.  
	\item If $N$ is a discrete-valued rv, then $Z_N$ is separable and we can choose a countable 
	separant $\diad$ for both $Z$ and $Z_N$ and further  
	$$ \pk*{\supA \norm{Z_N(t)}>0}= \pk*{ \norm{\ZT(0)}>0}=q.$$
If $N$ is a $\TT$-valued rv with positive pdf $\gamma(t)>0, t\in \TT$, then 
$Z_N$ does not need to be separable and moreover it is possible, if $Z$ is 
not stochastically continuous, that 
$$\pk*{\supA \norm{Z_N(t)}>0}=q <1 .$$ 
Since for all $t\in \TT$
$$\E{\norm{Z(t)}^\alpha}= \E{\norm{Z_N(t)}^\alpha}=1, 
$$
we have 
$q \in (0,1)$.   In view of \cite[Thm 9.4.2]{MR1280932} the joint measurability of  $Z_N\in \clasP$ 
  implies that it has a version (denote  it again by $Z_N$), 
  which is jointly measurable and separable   with a countable separant $\diad$ as $Z$.\\
  Define the RF 
 \bqn{\label{rL} 
  \widehat{ Z_N(t)} &=& q^{1/\alpha}Z_N(t) \Bigl \lvert \Bigl( 
	\sup_{s\in \TT} \norm{Z_N(s)}>0 \Bigr), \quad t\in \TT,
  }  
  which is jointly measurable and separable, since $Z_N$ is (as mentioned above). Since 
  $$ \pk*{\supA \norm{\widehat{Z_N(t)}}>0}=1,$$
   then \Cref{Theo2} implies that 
   $\widehat{ Z_N} \in \kk$. 

   \item The construction of $Z_N''$ is new also for the case of $\kk$ with 
 stochastically continuous elements satisfying the additional assumption 
 \eqref{specC}.
 \end{enumerate}
 \ERM

\section{Discussions and Extensions }
Under the assumptions and the notation of the previous section, in this  section we shall discuss the following:  
\begin{enumerate}[D1)]
	\item \label{iD1} alternative choices for  $\norm{\cdot}$;
	\item \label{iD2} generating shift-invariant classes via 
	the  tail and spectral tail RFs; 
	\item \label{iDtail} tail measures and \(L^\alpha\)-regular representors;
     \item \label{iD3} the special case of max-stable stationary RFs $X_Z$ with locally bounded sample paths;
	 \item \label{iD4} constructions of shift-generated classes via cluster RFs.\\   
\end{enumerate}

\underbar{\Cref{iD1} Choice of the norm}: Key properties
of the $\norm{\cdot}$ used in our definitions are  its 1-homogeneity 
and  measurability. As in \cite{hashorva2021shiftinvariant} 
we can utilise instead a general 1-homogeneous 
map $\norm{\cdot} \in \mathcal{H}_1 $ such that $\norm{B^h Z}, h\in \TT$ is jointly measurable.  
For such a general map we cannot show in general the existence of 
$Z^* \in \kk$ being $L^\alpha$-continuous. 
Since $\mathcal C [\norm{Z}]$ is also a shift-generated class
 of $1$-dimensional RFs and $S(\ZT)$ has the same
  law as $S(\norm{\ZT})$ we can still prove that  both \eqref{pandm} 
  and \eqref{e11} hold. \\

\underbar{\Cref{iD2} Tail and spectral tail RFs}:
\BD We shall call  $ \Theta \in \clasP$  
a spectral tail RF if
\begin{enumerate}[(i)]
	\item \label{defA:0} $\pk{\norm{\Theta(0)}=1}=1$;
	\item \label{defA:1} $\pk{\mathcal{S}(\Theta)>0}=1$;
	\item \label{defA:2} for all   $\Gamma\in  \Hh_0,  h\in \TT$
     the functional identity \eqref{eqDo20}. 
 \end{enumerate}    
\label{defA}
\ED

Based on our previous results, it is clear that any local RF
 $\widetilde \Theta$ of a given $\kk$ is a spectral tail RF. 
The converse also holds:   

\BT  
Any spectral tail RF $\Theta$ is a local RF of some $\kk$. 
Moreover, there exists an $L^\alpha$-continuous spectral tail RF $\Theta^\star  \in \clasP$ 
that has the same  finite dimensional distributions (fidis) as $\Theta$.  
\label{thT}
\ET
 
Recall that the positive rv $R$, independent of all other random elements, 
has survival function $\pk{R > s} = s^{-\alpha}$ for $s \ge 1$.   
Define the tail RF by
\[
Y(t) = R \Theta(t),
\]
a construction also referred to as the \emph{tail RF}; see, e.g., \cite{kulik:soulier:2020}.   
 
As first shown in \cite{Hrovje}, the following fundamental identity holds:
\bqn{
\E{ \Gamma( xB^h Y) \ind{ x \norm{ Y(-h)} > 1 } } 
&=& x^\alpha \ETI{ \Gamma( \widetilde{Y} ) \ind{ \norm{ \widetilde{Y}(h)} > x } }, 
\quad \forall \Gamma \in \mathcal{H},\ \forall h \in \TT,\ \forall x > 0.
\label{tYY}
}

 From our  findings above, it follows that identity \eqref{tYY} extends naturally to general maps 
of the form $F = F_1 F_2$, where $F_1 \in \mathcal{H}$ and $F_2 \in \Hh_\beta$ for some $\beta \ge 0$. 
In particular, we may take $F_2$ to be an integral functional.   
 
As discussed in \cite{hashorva2021shiftinvariant,MartinE}, the tail RF can also be defined directly as follows:

\BD We shall  call  a RF  $\vk Y \in \clasP   $ a   tail RF if
\begin{enumerate}[(i)]
	\item \label{defB:0} $\pk{\norm{ Y(0)}> 1}=1$;
	\item \label{defB:1} $\pk{\mathcal{S}(Y)> 0}=1$;
	\item \label{defB:2} \eqref{tYY} holds with $\widetilde{Y}$ substituted by $Y$.  
\end{enumerate}      
\label{defB}
\ED 
Interestingly, the functional equations that define a tail RF $Y$ can be extended as follows:

\begin{lemma}
	 If $Y$ is a tail RF, then   \eqref{tYY} holds   for $\widetilde{Y}=Y$ and all $x>0$ and all 
	 $F =F_1 F_2 F_3$ with $F_1 \in \mathcal{H}, F_2 
	 \in \Hh_\beta, \beta\ge 0, F_3(f)= (\mathcal{S}_{\mathbb{I}}(f))^a, a\inr, f\in \Dset$.
\label{kkS}
	\end{lemma}

 Tail and spectral tail RFs are crucial for the study of  stationary regularly varying time series, see e.g., 
 \cite{MR3561100, SegersEx17, klem,kulik:soulier:2020,MartinE, Resnickart,Ilya25}. 
In \cite{MartinE, Guenter} these RFs are introduced in an abstract setting. 

\BT 
If \(\Theta\) is the local RF of a shift-generated class \(\kk\), and
\(R\) is independent of \(\Theta\) with
\[
        \pk*{R>s}=s^{-\alpha},\qquad s\ge1,
\]
then
$
        Y=R\Theta
$
is a tail RF. Conversely, if \(Y\) is a tail RF, then
$
        \Theta=Y/ {\norm{Y(0)}}
$
is a spectral tail RF. In particular, there exists a shift-generated
class \(\kk\) such that \(Y/\norm{Y(0)}\) is its local RF.
\label{thmL}
\ET

\BRM   
	If $\Theta$ is a spectral tail RF, then 
      $$\pk{ \mathcal{S}(\Theta)>0 }=1 \quad \implies 
	  \pk{\mathcal{S}_{\mathbb{I}}(\Theta)>0}=1,$$ 
      since  we can find a stochastically continuous $\Theta^*$ 
      with the same fidis as $\Theta$, which is stochastically 
      continuous and the rest of the arguments are as in the proof 
      of \Cref{thmAug23}.
      On the other side, if in the definition of the spectral RF 
      we substitute $\pk{ \mathcal{S}(\Theta)>0 }=1$ with 
      $\pk{\mathcal{S}_{\mathbb{I}}(R \Theta)>0 }=1$,
       then since $\pk{\norm{\Theta(0)}=1}=1$, 
       it follows that $\pk{ \mathcal{S}(\Theta)>0 }=1$.
       The same reasoning applies when considering the definition of the tail RFs. 
   
\ERM

\underbar{\Cref{iDtail} Tail measures and \(L^\alpha\)-regular representors}:
Assume   as in Section 2 that
$
        \E{\|Z(0)\|^\alpha}=1 .$
For \(Z\in\kk\), define the   measure \(\nu_Z\) on
\((\Dset,\AA)\) by
\[
        \nu_Z(A)
        =
        \E*{
        \int_0^\infty
        \ind{rZ\in A}\,\alpha r^{-\alpha-1}\,dr
        },
        \qquad A\in\AA , 
\] 
which is \(-\alpha\)-homogeneous i.e.,  
\[
        \nu_Z(cA)=c^{-\alpha}\nu_Z(A),
        \qquad c>0,\quad A\in\AA .
\]

The measure \(\nu_Z\) depends only on the shift-generated class \(\kk\), and
not on the particular representor chosen in \(\kk\), see e.g., \cite{MartinE}, i.e.,  
\[
        \nu_Z=\nu_{\ZT},\qquad \ZT\in\kk 
\]
and  since \(B^hZ\in\kk\) it is shift-invariant, meaning that  
\[
        \nu_Z(H\circ B^h)
        =
        \nu_{B^hZ}(H)
        =
        \nu_Z(H),
        \qquad h\in\TT .
\] 

Our findings imply that  \(\nu_Z\) always admits
an \(L^\alpha\)-regular representor. Indeed, by \Cref{thmAug23}, there exists
\(Z^*\in\kk\) which is \(L^\alpha\)-continuous and hence  
\[
        \nu_Z=\nu_{Z^*}.
\]
Thus, although the initial representor \(Z\) is only assumed to be jointly
measurable and separable, the corresponding shift-invariant tail measure can always be
represented by an \(L^\alpha\)-continuous element of the same
shift-generated class.

The same regularization holds at the level of the local 
RFs. Let \(\Theta^*\) be the local RF associated with \(Z^*\). By
\Cref{thmAug23}, \(\Theta^*\) is \(L^\alpha\)-continuous. Setting 
\[
        Z_\gamma^*(t)
        =
        \frac{B^N\Theta^*(t)}
        {(\mathcal S_\gamma(B^N\Theta^*))^{1/\alpha}},
        \qquad t\in\TT 
\]
by \Cref{Theo2} and  \Cref{rLL}, \(Z_\gamma^*\in\kk\)  
for every non-negative measurable \(H\)
\[
\begin{aligned}
        \nu_Z(H)
        &=
        \E*{
        \int_0^\infty
        H(rZ_\gamma^*)\,\alpha r^{-\alpha-1}\,dr }                                                  \\
        &=
        \int_{\TT}
        \E*{
        \frac{1}{\mathcal S_\gamma(B^h\Theta^*)}
        \int_0^\infty
        H(rB^h\Theta^*)\,\alpha r^{-\alpha-1}\,dr
        }
        \gamma(h)\lambda(dh).
\end{aligned}
\]
Consequently, the tail measure \(\nu_Z\) is represented entirely in terms of the
\(L^\alpha\)-continuous spectral tail RF \(\Theta^*\). 
If, in addition, \(\|Z^*(0)\|>0\) almost surely, then one also has the
simpler   representation
\[
        \nu_Z(H)
        =
        \E*{
        \int_0^\infty
        H(r\Theta^*)\,\alpha r^{-\alpha-1}\,dr },
\] 
In general, when
\(\pk{\|Z^*(0)\|=0}>0\), this direct representation at the origin is not
available, and the randomized-shift representation above is the appropriate
substitute.

Consequently, \(L^\alpha\)-regularity is intrinsic to the tail measure
\(\nu_Z\).\\

\underbar{\Cref{iD3} $X_Z$ with locally bounded sample paths}: In \cite{hashorva2021shiftinvariant} it has been further assumed that  the max-stable stationary  $X_Z$ has locally bounded sample paths, i.e., $Z$ satisfies 
\bqn{\label{boundsp}
\E*{\sup_{t\in K} \norm{Z(t)}^\alpha} < \IF
} 
for all compact sets $K \subset \TT$. 
This assumption implies in particular that if for some $\ZT \in \kk$ 
\bqn{
	\pk{\SZT < \IF}=1,
\label{last}
} 
then also $\pk{\mathcal{S}(Z^*) < \IF}=1$ for some 
stochastically continuous $Z^*\in \kk$. 
Hence from \cite{hashorva2021shiftinvariant} we have that 
for some local RF $Y$ (recall $\mathcal{S}_{\mathbb{I}}(Y)=\int_{\TT}\ind{\norm{Y(t)}>1}\lambda(dt)$)
\bqn{
	\pk{\mathcal{S}_{\mathbb{I}}(Y) < \IF}=1.
\label{last2}
}

Moreover, the converse also holds; that is, under assumption \eqref{boundsp}, condition \eqref{last2} implies \eqref{last}.  

However, without assumption \eqref{boundsp}, we cannot conclude in general that \eqref{last} implies \eqref{last2}.  

Indeed, as shown in \cite{DombKabA}, there exists a process $Z$ such that \eqref{last} holds, 
yet assumption \eqref{boundsp} is not satisfied. \\

\underbar{\Cref{iD4} Constructions of shift-generated classes via cluster RFs}: 
Cluster RFs have been discussed recently  for stochastically continuous RFs in \cite{hashorva2025cluster}.  They have been shown to be instrumental for constructing  non-empty shift-generated classes $\kk$ such that \eqref{boundsp} holds and further $\pk{\mathcal{S}(Z) < \IF}=1$. When the latter condition is satisfied, taking $N$ as above, we have that $Z_N \in \kk$ 
with 
$$ Z_N(t) =   (\gamma(N))^{-1/\alpha} B^N Q(t), \quad t\in \TT ,$$
where $Q(t) = \Theta(t)/(\mathcal{S}(\Theta))^{1/\alpha}, t\in \TT$. The RF $Q$ is not a spectral tail RF, but it satisfies 
\bqn{ \label{q1}
	\pk{\mathcal{S}(Q)=1}=1. 
} 
\BD We call $Q \in \clasP$ a cluster RF, 
if it satisfies \eqref{q1} and $\pk*{\sup_{t\in \TT } \norm{Q(t)}>0}=1$.
 \ED  
 Given a cluster RF $Q$, we can define $Z_N$ as above. Recall that as noted in 
 Remark \ref{rLL} we can find a joint measurable and separable version of $Z_N$ with separant $\diad$. It can however be that 
 $\pk{\supA \norm{Z_N(t)}>0}< 1$. In that case we consider $\widehat{Z_N}$ defined in \eqref{rL} and then consider $\mathcal{C}[\widehat{Z_N}]$. One  interesting instance when $\alpha=1$ is by taking  $Q$ to be a continuous and strictly positive pdf  on $\TT$. In this case $\diad$ can be any countable dense subset of $\TT$.\\

When $\pk{\mathcal{S}(Z) < \IF}=1$,  it is possible to construct different cluster RFs for a given $\kk$, since they are not unique.  The advantage of this is that alternative formulas for extremal indices can be obtained, see \cite{kulik:soulier:2020, hashorva2025cluster}.  
Cluster random fields are also important in numerous statistical applications, see e.g., \cite{kulik:soulier:2020,MR4280158,chen2025asymptotic,chen2023limit}.

\section{Proofs}

	\prooftheo{verbpars}
	We first show that \eqref{boll} implies \eqref{boll22}. 
Fix \(h\in\TT\), \(\ZT\in\kk\), and \(G\in\Hh_0\). 
Let \(U_1=Z\) and \(U_2=B^h\ZT\). On a product enlargement define a mixture RF \(U\) by
\[
        U=
        \begin{cases}
        U_1, & \text{with probability }1/2,\\
        U_2, & \text{with probability }1/2.
        \end{cases}
\]
Applying the definition of \(\Hh_0\) to this mixture gives one sequence of
finite-dimensional \(0\)-homogeneous approximants. Thus, for suitable points
\(t_1^{(n)},\ldots,t_n^{(n)}\in\TT\) and Borel measurable
\(0\)-homogeneous maps \(G_n:(\R^d)^n\to[0,\IF]\),
\[
        G_n\bigl(U(t_1^{(n)}),\ldots,U(t_n^{(n)})\bigr)
        \toprob G(U), \qquad n\to\IF.
\]
Consequently, the same sequence approximates both \(G(Z)\) and
\(G(B^h\ZT)\) in probability.

For \(k>0\), put
\[
        G_{n,k}(f)
        =
        G_n\bigl(f(t_1^{(n)}),\ldots,f(t_n^{(n)})\bigr)\wedge k,
        \qquad f\in\Dset .
\]
Then \(G_{n,k}\in\mathcal H_0\), and as $n\to\IF$
\[
        G_{n,k}(Z)\toprob G(Z)\wedge k,
        \qquad
        G_{n,k}(B^h\ZT)\toprob G(B^h\ZT)\wedge k .
\]
Since
\[
        f\mapsto \|f(h)\|^\alpha G_{n,k}(f)
\]
belongs to \(\mathcal H_\alpha\), \eqref{boll} gives
\[
        \E{\|Z(h)\|^\alpha G_{n,k}(Z)}
        =
        \ETI{\|\ZT(0)\|^\alpha G_{n,k}(B^h\ZT)}.
\]
Letting \(n\to\IF\), and using boundedness of \(G_{n,k}\) together with
\[
        \E{\|Z(h)\|^\alpha}=1,
        \qquad
        \ETI{\|\ZT(0)\|^\alpha}=1,
\]
yields
\[
        \E{\|Z(h)\|^\alpha (G(Z)\wedge k)}
        =
        \ETI{\|\ZT(0)\|^\alpha (G(B^h\ZT)\wedge k)}.
\]
Finally, letting \(k\to\IF\) and using monotone convergence gives
\[
        \E{\|Z(h)\|^\alpha G(Z)}
        =
        \ETI{\|\ZT(0)\|^\alpha G(B^h\ZT)},
\]
which is \eqref{boll22}.\\

We next show that \eqref{boll22} implies \eqref{boll2}. 
Let \(F\in\Hh_\alpha\),   assume for simplicity that \(\diad\) is a common
separant of \(Z\) and \(\ZT\), and set
\[
        V(f)=\sum_{s\in\diad}q(s)\|f(s)\|^\alpha,
        \qquad q(s)>0,\quad \sum_{s\in\diad}q(s)=1.
\]
Then
\[
        V(Z),V(\ZT)\in(0,\IF)
        \quad\text{a.s.},
        \qquad
        \E {V(Z)}=\ETI {V(\ZT)}=1.
\]
With the convention \(F(f)/V(f)=0\) on \(\{V(f)=0\}\), put
\[
        G(f)=\frac{F(f)}{V(f)}.
\]
By the definition of \(\Hh_\alpha\) and the finite-dimensional approximation
of \(V\) through its partial sums, \(G\in\Hh_0\).

For each \(t\in\diad\), apply \eqref{boll22} with \(B^{-t}\ZT\) in place of
\(\ZT\). Since \(B^{-t}\ZT\in\kk\), we obtain
\[
        \E{\|Z(t)\|^\alpha G(Z)}
        =
        \ETI{\|\ZT(t)\|^\alpha G(\ZT)}.
\]
Therefore, by Tonelli's theorem,
\[
\begin{aligned}
        \E {F(Z)}
        &=
        \E{V(Z)G(Z)}  \\
        &=
        \sum_{t\in\diad}q(t)
        \E{\|Z(t)\|^\alpha G(Z)}  \\
        &=
        \sum_{t\in\diad}q(t)
        \ETI{\|\ZT(t)\|^\alpha G(\ZT)}  \\
        &=
        \ETI{V(\ZT)G(\ZT)}  \\
        &=
        \ETI {F(\ZT)}.
\end{aligned}
\]
Since \(B^h\ZT\in\kk\) for every \(h\in\TT\), replacing \(\ZT\) by
\(B^h\ZT\) gives
\[
        \E {F(Z)}=\ETI {F(B^h\ZT)},
        \qquad h\in\TT,
\]
which is \eqref{boll2}. The implication \eqref{boll2}\(\Rightarrow\)\eqref{boll}
is immediate.
\\ 	
It remains to justify the claims involving the integral functional. Suppose for simplicity that \(\gamma(t)\equiv 1\). 

For \(U\in\clasP\),  by \Cref{lemGen2} applied to
\(W(t)=\|U(t)\|^\alpha\)   there exist
finite-dimensional \(\alpha\)-homogeneous maps \(A_k\) such that, after
passing to a subsequence,
\[
        A_k(U)\to \mathcal S (U)
        \quad\text{a.s.}
\]

For \(a=0\), the sampling points in \Cref{lemGen2} can be chosen so that
\[
        \ind{\{A_k(U)=0\}}
        \to
        \ind{\{\mathcal S(U)=0\}}
        \quad\text{a.s.}
\]
Indeed, with \(A_0=\{\mathcal S(U)=0\}\), Tonelli's theorem gives
\[
        \int_{\TT}
        \mathbb P\{A_0,\|U(t)\|>0\}\,\lambda(dt)=0.
\]
Hence the deterministic sampling points may be chosen outside the exceptional
\(\lambda\)-null set.

For \(a=\IF\), ordinary thresholds such as \(\ind{\{A_k>r\}}\) are not
\(0\)-homogeneous. We therefore use a scale-free ratio. Define
\[
        C_k(f)=\max_{1\le i\le k}A_i(f),
        \qquad
        D_k(f)=\sum_{i=1}^k2^{-i}C_i(f),
\]
and
\[
        R_k(f)=\frac{C_k(f)}{D_k(f)},
        \qquad 0/0:=0.
\]
Then \(R_k\) is finite-dimensional and \(0\)-homogeneous. Moreover,
\[
        R_k(U)\to\IF
        \quad\text{on}\quad
        \{\mathcal S(U)=\IF\}, \qquad k \to\infty 
\]
whereas \(R_k(U)\) remains finite on
\(\{\mathcal S(U)<\IF\}\). Hence one can choose a subsequence \(k_n\)
such that
\[
        H_n(f)=\ind{\{R_{k_n}(f)>n\}}
\]
satisfies as \(n\to\IF\)
\[
        H_n(U)\toprob
        \ind{\{\mathcal S(U)=\IF\}}.
\]
The maps \(H_n\) are finite-dimensional and \(0\)-homogeneous.

Multiplying these \(0\)-homogeneous approximants by finite-dimensional
\(\beta\)-homogeneous approximants of \(\Gamma_\beta\) gives
\[
        \Gamma_\beta(f)\ind{\{\mathcal S(f)=a\}}
        \in\Hh_\beta,
        \qquad a\in\{0,\IF\}.
\]
Similarly, using the approximants \(A_k\) and a diagonal argument gives
\[
        \Gamma_\beta(f)\mathcal S(f)\in\Hh_{\alpha+\beta}.
\]
Finally, with the convention
\[
        \frac{\Gamma_\alpha(f)}{\mathcal S(f)}=0
        \quad\text{on}\quad
        \{\mathcal S(f)\in\{0,\IF\}\},
\]
the same approximation and diagonal argument gives
$
        \Gamma_\alpha(f)/\mathcal S(f)\in\Hh_0.
$ 
	\QED
	\\

\prooftheo{thmAug23}
We consider for simplicity \(d=1\), the general case follows with similar arguments. The vector-valued case follows by
applying the same construction component-wise.

By \Cref{verbpars}, since \(Z\in\kk\), we have
\[
        \E{\norm{Z(h)}^\alpha G(Z)}
        =
        \E{\norm{Z(0)}^\alpha G(B^h Z)},
        \qquad h\in\TT,\quad G\in\Hh_0 .
\]

Set
\[
        Z_+(t)=\max(Z(t),0),
        \qquad
        Z_-(t)=\max(-Z(t),0),
\qquad        \bar Z(t)=(Z_+(t),Z_-(t)),\qquad t\in\TT .
\]
and equip \(\mathbb R^2\) with the norm
$
        \|(x_1,x_2)\|_\infty=\max(|x_1|,|x_2|), $ hence 
$ 
        \norm{\bar Z(t)}_\infty=|Z(t)|.$

Let \(X\) be the two-dimensional max-stable RF generated by \(\bar Z\),
namely
\[
        X(t)=
        \max_{i\ge1}
        \frac{\bar Z^{(i)}(t)}
        {\left(\sum_{k=1}^i\mathcal V_k\right)^{1/\alpha}},
        \qquad t\in\TT,
\]
where the maximum is component-wise. The RF \(X\) is well-defined since
\[
        \E{Z_+(t)^\alpha+Z_-(t)^\alpha}<\infty,
        \qquad t\in\TT .
\]

We show that \(X\) is stationary. Let \(F\in\mathcal H_0\), and define
\(G_0(z)=(z_+,z_-)\). Since \(G_0\) is measurable and \(1\)-homogeneous,
\(F\circ G_0\in\mathcal H_0\). Hence
\[
\begin{aligned}
\E{\norm{\bar Z(h)}_\infty^\alpha F(\bar Z)}
&=
\E{|Z(h)|^\alpha F(G_0(Z))}                                      \\
&=
\E{|Z(0)|^\alpha F(G_0(B^h Z))}                                  \\
&=
\E{\norm{\bar Z(0)}_\infty^\alpha F(B^h\bar Z)} .
\end{aligned}
\]
By the   tilt characterization of stationarity for max-stable RFs established in \cite{Htilt}
\(X\) is stationary. Moreover, the maximum of jointly measurable and separable RFs $\bar Z^{(i)}$ is jointly measurable and separable. 
In view of \cite[Thm 3.1]{StoevSPA} or  \cite[Thm 1.3.3]{MR3561100},  
this is equivalent with $X $ having a version $\widetilde{X}$, 
which is  stochastically continuous. Note in passing that the latter theorem assumes finiteness of the second moments of $X$. Since $X'=X^\beta$ is also max-stable for $\beta>0$, we can apply the theorem to $X'$ with an appropriate $\beta$. \\ 
 
By \cite[p.~1765]{dom2016}, \(\widetilde X\) has an
\(L^\alpha\)-continuous non-negative representor
$ 
        \widetilde{\bar Z}=(\ZT_1,\ZT_2)$ and hence \(\ZT_1\) and \(\ZT_2\) are \(L^\alpha\)-continuous.  Remark  that the aforementioned result is for $d=1$ and the extension to $d>1$ follows by similar arguments. \\

Set
\[
        \Delta=\ZT_1-\ZT_2 
\]
and let  \(\diad\) be a common countable separant for \(Z\) and \(\Delta\),
enlarged if necessary so that \(0\in\diad\). Define
\[
        A=
        \left\{
        \sup_{h\in\diad}|\Delta(h)|>0
        \right\},
        \qquad
        p=\pk{A}.
\]
 By the representor-equivalence, e.g.,  \cite[Prop 2.1]{KumeE} applied to
the \(\alpha\)-homogeneous functional
\[
        H(u,v)=|u(0)-v(0)|^\alpha ,
\]
we have
\[
        \E{|\Delta(0)|^\alpha}
        =
        \E{|Z_+(0)-Z_-(0)|^\alpha}
        =
        \E{|Z(0)|^\alpha}
        =
        1.
\]
Hence \(\pk{|\Delta(0)|>0}>0\). Since \(0\in\diad\), this implies \(p>0\).

If  \(Z^*\) denotes  the conditional RF 
$
        p^{1/\alpha}\Delta\lvert  A,$  then
$
        \pk*{\sup_{h\in\diad}|Z^*(h)|>0}=1 .
$ We prove that \(Z^*\in\kk\). Let \(F\in\mathcal H_\alpha\) and \(h\in\TT\).
Define
\[
        H_{F,h}(u,v)
        =
        \begin{cases}
        F(B^h(u-v)),&
        \displaystyle \sup_{r\in\diad}|u(r)-v(r)|>0,\\[1mm]
        0,&
        \displaystyle \sup_{r\in\diad}|u(r)-v(r)|=0.
        \end{cases}
\]
Then \(H_{F,h}\) is measurable and \(\alpha\)-homogeneous. The second line
only fixes the value on the zero-reconstruction event.

By the definition of \(Z^*\) and the \(\alpha\)-homogeneity of \(F\) we have 
\[
\begin{aligned}
        \E{F(B^h Z^*)}
&=
        \E*{F\left(B^h(p^{1/\alpha}\Delta)\right)
        \lvert A}                                                   \\
&=
        p\E{F(B^h\Delta)\mid A}                            \\
&=
        \E{F(B^h\Delta)\ind{A}}                              \\
&= \E{H_{F,h}(\ZT_1,\ZT_2)} .
\end{aligned}
\]
Since \((\ZT_1,\ZT_2)\) and \((Z_+,Z_-)\) are representors of the same
two-dimensional max-stable RF, the representor-equivalence theorem gives
\[
        \E{H_{F,h}(\ZT_1,\ZT_2)}
        =
        \E{H_{F,h}(Z_+,Z_-)} .
\]
But \(Z=Z_+-Z_-\), and by separability and \eqref{eq1} 
\[
        \pk*{\sup_{r\in\diad}|Z(r)|>0}=1
\] implying that almost surely $ 
        H_{F,h}(Z_+,Z_-)=F(B^h Z)$. 
Consequently, we have 
\[
        \E{F(B^h Z^*)}
        =
        \E{F(B^h Z)}
        =
        \E{F(Z)},
        \qquad h\in\TT,\quad F\in\mathcal H_\alpha .
\]
Thus \(Z^*\in\kk\).

We next show that \(Z^*\) is \(L^\alpha\)-continuous. For \(s,t\in\TT\),
\[
\begin{aligned}
\E{|Z^*(s)-Z^*(t)|^\alpha}
&=
p\,\E{|\Delta(s)-\Delta(t)|^\alpha\mid A}                   \\
&=
\E{|\Delta(s)-\Delta(t)|^\alpha\ind{A}}                     \\
&\le
\E{|\Delta(s)-\Delta(t)|^\alpha}.
\end{aligned}
\]
Since \(\Delta=\ZT_1-\ZT_2\),
\[
        |\Delta(s)-\Delta(t)|
        \le
        |\ZT_1(s)-\ZT_1(t)|+|\ZT_2(s)-\ZT_2(t)|
\]
hence, with
\[
        c_\alpha=
        \begin{cases}
        1, & 0<\alpha\le1,\\
        2^{\alpha-1}, & \alpha>1,
        \end{cases}
\]
we have
\[
\begin{aligned}
\E{|\Delta(s)-\Delta(t)|^\alpha}
&\le
c_\alpha
\left[
        \E{|\ZT_1(s)-\ZT_1(t)|^\alpha}
        +
        \E{|\ZT_2(s)-\ZT_2(t)|^\alpha}
\right].
\end{aligned}
\]
The right-hand side tends to \(0\) as \(s\to t\), because
\(\ZT_1,\ZT_2\) are \(L^\alpha\)-continuous. Hence \(Z^*\) is
\(L^\alpha\)-continuous. Taking a jointly measurable and separable version if
necessary, we still denote it by \(Z^*\).

Since \(Z^*\) is stochastically continuous and
\[
        \pk*{\sup_{h\in\diad}|Z^*(h)|>0}=1,
\]
\cite[Thm~9.1]{hashorva2021shiftinvariant}, or equivalently
\cite[Thm~2.1]{GeoSS}, yields
\[
        \pk{\mathcal S(Z^*)>0}=1 .
\]

We now transfer the positivity to any \(\ZT\in\kk\). Let \(\diad\) be a
common countable separant of \(Z^*\) and \(\ZT\), and choose
\(q(h)>0\), \(h\in\diad\), with
\[
        \sum_{h\in\diad}q(h)=1.
\]
Setting $ 
        V(f)=\sum_{h\in\diad}q(h)\|f(h)\|^\alpha .
$ we have  almost surely
$
        V(Z^*),\,V(\ZT)\in(0,\infty)$.  
Indeed, finiteness follows from Tonelli's theorem and the normalization, while
strict positivity follows from separability and \eqref{eq1}.

By \Cref{verbpars}, we have
\[
        f\mapsto V(f)\ind{\mathcal S(f)=0}
        \in\Hh_\alpha 
\] 
and hence 
\[
\begin{aligned}
0
&=
        \E{
        V(Z^*)\ind{\mathcal S(Z^*)=0}
        }                                                   =
        \ETI{
        V(\ZT)\ind{\mathcal S(\ZT)=0}
        } .
\end{aligned}
\]
Since \(V(\ZT)>0\) almost surely, then 
$
        \pTI{\mathcal S(\ZT)=0}=0$ and hence  
\[
        \pTI{\mathcal S(\ZT)>0}=1,
        \qquad \ZT\in\kk,
\]
which proves \eqref{pandm2}.

Let \(\Theta^*\) be the local RF associated with \(Z^*\). Under the tilted
probability associated with \(Z^*\),
\[
        \Theta^*(t)=\frac{Z^*(t)}{\|Z^*(0)\|}
\]
on \(\{\|Z^*(0)\|>0\}\), and we define it arbitrarily on the complementary
tilted-null set. Then, for \(s,t\in\TT\),
\[
\begin{aligned}
\E{\|\Theta^*(s)-\Theta^*(t)\|^\alpha}
&=
\E{
        \|Z^*(0)\|^\alpha
        \left\|
        \frac{Z^*(s)}{\|Z^*(0)\|}
        -
        \frac{Z^*(t)}{\|Z^*(0)\|}
        \right\|^\alpha
        \ind{\|Z^*(0)\|>0}
        }                                                   \\
&=
\E{
        \|Z^*(s)-Z^*(t)\|^\alpha
        \ind{\|Z^*(0)\|>0}
        }                                                   \\
&\le
\E{\|Z^*(s)-Z^*(t)\|^\alpha}.
\end{aligned}
\]
Here the expectation on the left is under the tilted law, whereas the last
expectation is under the original law of \(Z^*\). Hence \(\Theta^*\) is
\(L^\alpha\)-continuous.

Again by \cite[Thm~9.1]{hashorva2021shiftinvariant}, or
\cite[Thm~2.1]{GeoSS},
\[
        \pk{\mathcal S(\Theta^*)>0}=1 .
\]

Let now \(\widetilde\Theta\) be the local RF associated with an arbitrary
\(\ZT\in\kk\). For every bounded finite-dimensional \(0\)-homogeneous
functional \(\Gamma\),
\[
\begin{aligned}
\E{\Gamma(\widetilde\Theta)}
&=
        \ETI*{
        \|\ZT(0)\|^\alpha
        \Gamma\left(\frac{\ZT}{\|\ZT(0)\|}\right)
        }                                                   \\
&=
        \ETI{
        \|\ZT(0)\|^\alpha \Gamma(\ZT)
        }                                                   \\
&=
        \E{
        \|Z^*(0)\|^\alpha \Gamma(Z^*)
        }                                                   \\
&=
        \E{\Gamma(\Theta^*)}.
\end{aligned}
\]
Thus all local RFs of \(\kk\) have the same finite-dimensional distributions
as \(\Theta^*\). By \Cref{lemGen2}, the law of
\(\mathcal S(\widetilde\Theta)\) is determined by these finite-dimensional
distributions implying 
\[
        \pk{\mathcal S(\widetilde\Theta)>0}
        =
        \pk{\mathcal S(\Theta^*)>0}
        =
        1.
\]

Finally, let
\[
        Y^*(t)=R\Theta^*(t),\qquad t\in\TT,
\]
where \(R\) is independent of \(\Theta^*\) and satisfies
\[
        \pk{R>s}=s^{-\alpha},\qquad s\ge1.
\]
Since \(\Theta^*\) is stochastically continuous, so is \(Y^*\). Moreover,
\[
        \norm{Y^*(0)}=R>1
        \quad\text{a.s.}
\]
Therefore, by \cite[Thm~2.1]{GeoSS},
\[
        \pk{\mathcal S_{\mathbb I}(Y^*)>0}=1.
\]
If \(\widetilde Y=R\widetilde\Theta\), with \(R\) independent of
\(\widetilde\Theta\), then
\[
        \widetilde Y\stackrel{\mathrm{fidi}}{=}Y^*.
\]
By \Cref{lemGen2}, the law of
\(\mathcal S_{\mathbb I}(\widetilde Y)\) is determined by the finite-dimensional
distributions of \(\widetilde Y\). Consequently,
\[
        \pk{\mathcal S_{\mathbb I}(\widetilde Y)>0}
        =
        \pk{\mathcal S_{\mathbb I}(Y^*)>0}
        =
        1.
\]
Together with the positivity of \(\mathcal S(\widetilde\Theta)\), this proves
\eqref{e11}.

\QED
\\   
\proofkorr{Theo2}
 We note first that $\mathcal{S}_\gamma( B^N \widetilde{\Theta} )$
  is finite a.s.\ follows easily since \eqref{defA} implies that 
$$ \E{\norm{\Theta(h)}^\alpha}= \pk{\norm{ \Theta(-h)}\not=0}\le 1$$
and moreover the df of $\mathcal{S}_\gamma( B^N \widetilde{\Theta} )$
 does not depend on the choice of $\widetilde{\Theta}$. Hence using further 
 \eqref{e11}
$$\pk{\mathcal{S}_\gamma( B^N \widetilde{\Theta} )>0 }=
\pk{\mathcal{S}_\gamma( B^N {\Theta^*} )>0 }=1,$$
with $\Theta^*$ a stochastically continuous local RF. \\
Similarly,  the finiteness of $\mathcal{S}_\gamma( B^N \ZT )$ 
follows from $ \E{\norm{\ZT(h)}^\alpha}=1, h\in \TT$. Since 
$$
Z_h= B^h \ZT \in\kk, \quad  \forall h\in \TT,
$$
 from \eqref{pandm} and the positivity of $\gamma$ 
we have that 
$\mathcal{S}_\gamma( Z_h )>0$ a.s.\ implying that 
$\mathcal{S}_\gamma( B^N \ZT )>0$ a.s.\ since $N$ is independent of $\ZT$. \\
  Next, the proof follows from \Cref{thmAug23} and \Cref{verbpars}. 
  We give therefore only the details for $Z_N'$.   It suffices to show that the RF $Z_N'$   satisfies \eqref{boll2} for all $h\in \TT, F \in \mathcal{H}_0$.\\
Since   $\pk{ \norm{\Theta(0)} =1}=1$, applying
 \eqref{eqDo20}  and  the Fubin-Tonelli theorem, we obtain 
\bqn{\label{verst}
	\E{\norm{Z'_N(h)}^\alpha F(Z_N')} &=&\int_{\TT}  \E*{ \norm{B^t \Theta(h)}^\alpha \frac{ \norm{\Theta(0)}^\alpha F\left(  B^{t} {\Theta}	\right)}{ \mathcal{S}_\gamma(B^t{\Theta} )} } \gamma(t)\lambda(dt)\nonumber \\
	&=& 
\int_{\TT}  \E*{  \frac{ \norm{\Theta(t-h)}^\alpha F\left(  B^h {\Theta}
	\right)}{ \mathcal{S}_\gamma(B^h{\Theta} )} } \gamma(t)\lambda(dt)\nonumber \\
	&=&
	 \E*{F\left(  B^h {\Theta}\right)  \int_{\TT}  \frac{ \norm{\Theta(t-h)}^\alpha 	}{ \mathcal{S}_\gamma({B^h\Theta} )} \gamma(t)\lambda(dt) } \nonumber \\
 &=&
	 \E*{F\left(  B^h {\Theta}\right)   }
	  \nonumber \\
	  &=&
	 \E{\norm{Z(0)}^\alpha F(B^h Z)} 
}
establishing the claim.   
\QED 
\\ 

\prooftheo{thT}
Let first \(\Theta\) be a local RF of a shift-generated class \(\kk\).
Then
\[
        \pk*{\norm{\Theta(0)}=1}=1
\]
by the definition of the local RF. Moreover,
\[
        \pk*{\mathcal S(\Theta)>0}=1
\]
follows from \Cref{thmAug23}, and the identity \eqref{eqDo20} follows
from \Cref{verbpars}. Hence every local RF is a spectral tail RF.

Conversely, let \(\Theta\) be a spectral tail RF. Let \(N\) be independent
of \(\Theta\), with strictly positive pdf \(\gamma\) with respect to
Lebesgue measure on \(\TT\). Define
\[
        \mathcal A_\gamma(f)
        =
        \int_{\TT}
        \gamma(s)\ind{\norm{f(s)}>0}\lambda(ds).
\]
By \Cref{lemGen2}, applied to
\[
        W_U(s)=\gamma(s)\ind{\norm{U(s)}>0},
\]
we have \(\mathcal A_\gamma\in\Hh_0\). Hence, for every fixed
\(a\in\TT\),
\[
        f\mapsto \mathcal A_\gamma(B^a f)
\]
also belongs to \(\Hh_0\).

Since
$
        \pk*{\mathcal S(\Theta)>0}=1, $
we have
\[
        \lambda\{s\in\TT:\norm{\Theta(s)}>0\}>0
        \quad \text{a.s.}
\]
Thus, because \(\gamma>0\)
\[
        \pk*{\mathcal A_\gamma(B^h\Theta)>0}=1,
        \qquad h\in\TT
\]
and in particular
$
        \pk*{\mathcal A_\gamma(B^N\Theta)>0}=1.$

Define
\[
        Z^0(t)
        =
        \frac{B^N\Theta(t)}
        {\mathcal A_\gamma(B^N\Theta)^{1/\alpha}},
        \qquad t\in\TT.
\]
The RF \(Z^0\) is jointly measurable.

We prove first the weighted shift identity for bounded \(G\in\Hh_0\).
For fixed \(n\in\TT\), the maps
\[
        f\mapsto
        \frac{G(B^n f)}
        {\mathcal A_\gamma(B^n f)}
        \ind{\mathcal A_\gamma(B^n f)>0}
\]
and
\[
        f\mapsto
        \frac{G(B^{h+n} f)}
        {\mathcal A_\gamma(B^n f)}
        \ind{\mathcal A_\gamma(B^n f)>0}
\]
belong to \(\Hh_0\), by the closure of \(\Hh_0\) under shifts, products
and division by a non-negative \(0\)-homogeneous functional on its
positivity set.

Using \eqref{eqDo20} with shift \(h-n\), and then Tonelli's theorem,
we obtain
\[
\begin{aligned}
\E*{\norm{Z^0(h)}^\alpha G(Z^0)}
&=
\int_{\TT}\gamma(n)
\E*{
\frac{
\norm{\Theta(h-n)}^\alpha G(B^n\Theta)
}{
\mathcal A_\gamma(B^n\Theta)
}
}\lambda(dn)                                                \\
&=
\int_{\TT}\gamma(n)
\E*{
\frac{
\ind{\norm{\Theta(n-h)}>0}G(B^h\Theta)
}{
\mathcal A_\gamma(B^h\Theta)
}
}\lambda(dn)                                                \\
&=
\E*{
G(B^h\Theta)
\frac{
\int_{\TT}\gamma(n)\ind{\norm{\Theta(n-h)}>0}\lambda(dn)
}{
\mathcal A_\gamma(B^h\Theta)
}
}                                                           \\
&=
\E*{G(B^h\Theta)}
\end{aligned}
\] 
implying 
\[
        \E*{\norm{Z^0(h)}^\alpha G(Z^0)}
        =
        \E*{\norm{Z^0(0)}^\alpha G(B^hZ^0)}
\]
for every bounded \(G\in\Hh_0\) and every \(h\in\TT\).

For general \(G\in\Hh_0\), apply the preceding identity to
\(G_k=G\wedge k\), \(k\ge1\), and let \(k\to\infty\). By monotone
convergence
\[
        \E*{\norm{Z^0(h)}^\alpha G(Z^0)}
        =
        \E*{\norm{Z^0(0)}^\alpha G(B^hZ^0)},
        \qquad G\in\Hh_0,\ h\in\TT.
        \tag{1}
\]
Taking \(G\equiv1\) and \(h=0\), we get
\[
        \E*{\norm{Z^0(0)}^\alpha}=1.
\]

Moreover almost surely 
\[
        \mathcal S(Z^0)
        =
        \frac{\mathcal S(B^N\Theta)}
        {\mathcal A_\gamma(B^N\Theta)}
        =
        \frac{\mathcal S(\Theta)}
        {\mathcal A_\gamma(B^N\Theta)}
        >0 
\]
and hence  \(Z^0\) is non-zero as a RF  almost surely.

Let \(Z\) be a jointly measurable separable version of \(Z^0\). Since the
identities above are determined by finite-dimensional distributions, \(Z\)
also satisfies
\[
        \E*{\norm{Z(h)}^\alpha G(Z)}
        =
        \E*{\norm{Z(0)}^\alpha G(B^hZ)},
        \qquad G\in\Hh_0,\ h\in\TT,
\]
and
\[
        \E*{\norm{Z(0)}^\alpha}=1,
        \qquad
        \pk*{\mathcal S(Z)>0}=1.
\]
In particular, we have 
\[
        \pk*{\sup_{t\in\TT}\norm{Z(t)}>0}=1.
\]
Hence \(Z\in\clasP\) satisfies \eqref{eq1}. By \Cref{verbpars}, identity
\((1)\) is equivalent to the defining functional identity \eqref{boll}.
Therefore \(\mathcal C[Z]\) is a non-empty shift-generated class.

It remains to identify its local RF. Let \(\Theta_Z\) be the local RF
obtained from \(Z\). Since \(Z\) and \(Z^0\) have the same
finite-dimensional distributions, it is enough to compute with \(Z^0\).
Let \(\Gamma\in\Hh_0\) be bounded. Then, by the definition of the local RF,
the \(0\)-homogeneity of \(\Gamma\), and \eqref{eqDo20} with shift \(-n\),
\[
\begin{aligned}
\E*{\Gamma(\Theta_Z)}
&=
\E*{
\norm{Z^0(0)}^\alpha
\Gamma\left(\frac{Z^0}{\norm{Z^0(0)}}\right)
}                                                           \\
&=
\int_{\TT}\gamma(n)
\E*{
\frac{
\norm{\Theta(-n)}^\alpha\Gamma(B^n\Theta)
}{
\mathcal A_\gamma(B^n\Theta)
}
}\lambda(dn)                                                \\
&=
\int_{\TT}\gamma(n)
\E*{
\frac{
\ind{\norm{\Theta(n)}>0}\Gamma(\Theta)
}{
\mathcal A_\gamma(\Theta)
}
}\lambda(dn)                                                \\
&=
\E*{
\Gamma(\Theta)
\frac{
\int_{\TT}\gamma(n)\ind{\norm{\Theta(n)}>0}\lambda(dn)
}{
\mathcal A_\gamma(\Theta)
}
}                                                           \\
&=
\E*{\Gamma(\Theta)}.
\end{aligned}
\]
By monotone convergence, the same identity holds for all non-negative
\(\Gamma\in\Hh_0\).

Since further 
\[
        \pk*{\norm{\Theta_Z(0)}=1}
        =
        \pk*{\norm{\Theta(0)}=1}
        =
        1,
\]
this equality for bounded \(0\)-homogeneous finite-dimensional test
functionals implies equality of finite-dimensional distributions. Hence
\[
        \Theta_Z\stackrel{\mathrm{fidi}}{=}\Theta.
\]
Thus \(\Theta\) is a local RF of the shift-generated class
\(\mathcal C[Z]\).

Finally, applying \Cref{thmAug23} to \(\mathcal C[Z]\), there exists an
\(L^\alpha\)-continuous local RF \(\Theta^\star\in\clasP\). All local RFs
of the same shift-generated class have the same finite-dimensional
distributions. Therefore
\[
        \Theta^\star\stackrel{\mathrm{fidi}}{=}\Theta.
\]
Since every local RF is a spectral tail RF, \(\Theta^\star\) is an
\(L^\alpha\)-continuous spectral tail RF.
\QED
 \\

\prooflem{kkS}
The claim follows by  approximating $F_2(Y)$ and $F_3(Y)$ 
utilising \Cref{lemGen2}.  
\QED 
\\

 \prooftheo{thmL}
Assume first that \(\Theta\) is the local RF of a shift-generated class
\(\kk\). By \Cref{thT}, \(\Theta\) is a spectral tail RF. Let \(R\) be
independent of \(\Theta\), with
$
        \pk*{R>s}=s^{-\alpha}, s\ge1,
$
and put
$
        Y=R\Theta .
$
Then almost surely 
$
        \norm{Y(0)}=R>1$   and 
$
        \mathcal S(Y)=R^\alpha \mathcal S(\Theta)>0$. 
It remains to verify \eqref{tYY}.

We first note that in view of \eqref{eqDo20},  for every
non-negative measurable functional \(H\), every \(h\in\TT\), and with the
usual convention on \(\{\norm{f(h)}=0\}\),
\bqn{\label{Pn}
        \E*{
        \norm{\Theta(h)}^\alpha
        H\left(\frac{\Theta}{\norm{\Theta(h)}}\right)}
        =
        \E*{
        \ind{\norm{\Theta(-h)}\ne0}H(B^h\Theta)} .}

Let now \(\Gamma\in\mathcal H\), \(h\in\TT\), and \(x>0\). Conditioning
on \(R\) and using the change of variables \(u=xr\), we get
\[
\begin{aligned}
\E*{\Gamma(xB^hY)\ind{x\norm{Y(-h)}>1}}
&=
\int_1^\infty
\alpha r^{-\alpha-1}
\E*{
\Gamma(xrB^h\Theta)
\ind{xr\norm{\Theta(-h)}>1}
}\,dr                                                     \\
&=
x^\alpha
\E*{
\int_x^\infty
\alpha u^{-\alpha-1}
\Gamma(uB^h\Theta)
\ind{u\norm{\Theta(-h)}>1}
\,du
}.
\end{aligned}
\]
Define
\[
        H_{x,\Gamma}(f)
        =
        \int_x^\infty
        \alpha u^{-\alpha-1}
        \Gamma(uf)\ind{u\norm{f(0)}>1}\,du .
\]
Then the last display equals
\[
        x^\alpha\E*{H_{x,\Gamma}(B^h\Theta)}.
\]
Using \eqref{Pn}, this is
\[
        x^\alpha
        \E*{
        \norm{\Theta(h)}^\alpha
        H_{x,\Gamma}\left(\frac{\Theta}{\norm{\Theta(h)}}\right)
        }.
\]
For \(a=\norm{\Theta(h)}\), a change of variables \(u=ar\) yields
\[
\begin{aligned}
a^\alpha
H_{x,\Gamma}\left(\frac{\Theta}{a}\right)
&=
\int_1^\infty
\alpha r^{-\alpha-1}
\Gamma(r\Theta)\ind{ra>x}\,dr .
\end{aligned}
\]
Consequently, we have 
\[
\begin{aligned}
\E*{\Gamma(xB^hY)\ind{x\norm{Y(-h)}>1}}
&=
x^\alpha
\E*{
\int_1^\infty
\alpha r^{-\alpha-1}
\Gamma(r\Theta)\ind{r\norm{\Theta(h)}>x}\,dr
}                                                        \\
&=
x^\alpha
\E*{\Gamma(Y)\ind{\norm{Y(h)}>x}} .
\end{aligned}
\]
Thus \(Y\) satisfies \eqref{tYY}, and hence \(Y\) is a tail RF.
 Conversely, let \(Y\) be a tail RF and put
$
        R=\norm{Y(0)},
        \Theta=\frac{Y}{\norm{Y(0)}} .
$
Taking \(h=0\) and \(\Gamma\equiv1\) in \eqref{tYY}, for \(x\ge1\), gives
\[
        1
        =
        \pk*{x\norm{Y(0)}>1}
        =
        x^\alpha\pk*{\norm{Y(0)}>x}.
\]
Hence
\[
        \pk*{R>x}=x^{-\alpha},
        \qquad x\ge1.
\]
Next, for any non-negative measurable \(F\), apply \eqref{tYY} with
\(h=0\) and
\[
        \Gamma(f)=F\left(\frac{f}{\norm{f(0)}}\right).
\]
For \(x\ge1\), we obtain
\[
        \E*{F(\Theta)}
        =
        x^\alpha\E*{F(\Theta)\ind{R>x}} .
\]
Thus \(R\) is independent of \(\Theta\)
 and almost surely
$
        \norm{\Theta(0)}=1$  
and $ 
        \mathcal S(\Theta)
        =
        R^{-\alpha}\mathcal S(Y)>0$.  
It remains to prove \eqref{eqDo20}. Given  \(\Gamma\in\Hh_0\),   \eqref{tYY} gives
\[
        \E*{
        \Gamma(B^h\Theta)
        \ind{xR\norm{\Theta(-h)}>1}
        }
        =
        x^\alpha
        \E*{
        \Gamma(\Theta)
        \ind{R\norm{\Theta(h)}>x}
        } .
\]
Letting \(x\to\infty\), the left-hand side converges by monotone
convergence to
\[
        \E*{
        \ind{\norm{\Theta(-h)}\ne0}
        \Gamma(B^h\Theta)
        }.
\]
For the right-hand side, using the independence of \(R\) and \(\Theta\),
\[
\begin{aligned}
x^\alpha
\E*{
\Gamma(\Theta)
\ind{R\norm{\Theta(h)}>x}
}
&=
\E*{
\Gamma(\Theta)
x^\alpha
\pk*{R>x/\norm{\Theta(h)}\mid\Theta}
}                                                       \\
&=
\E*{
\Gamma(\Theta)
\left(x^\alpha\wedge\norm{\Theta(h)}^\alpha\right)
}.
\end{aligned}
\]
Again by monotone convergence, this tends to
\[
        \E*{
        \norm{\Theta(h)}^\alpha\Gamma(\Theta)
        }.
\]
Therefore
\[
        \E*{
        \norm{\Theta(h)}^\alpha\Gamma(\Theta)
        }
        =
        \E*{
        \ind{\norm{\Theta(-h)}\ne0}
        \Gamma(B^h\Theta)
        },
        \qquad
        \Gamma\in\Hh_0,\ h\in\TT.
\]
Hence \(\Theta\) is a spectral tail RF. By \Cref{thT}, there exists a
shift-generated class \(\kk\) such that \(\Theta=Y/\norm{Y(0)}\) is its
local RF.
\QED
\\

\section{Appendix}
\BEL
Let \(W(t),t\in\TT=\R^l\), be a non-negative jointly measurable RF
defined on a complete probability space
\((\Omega_1,\mathscr F_1,\mathbb P_1)\). Define, for Borel sets
\(A\subset\TT\),
\[
        I_A(W)=\int_A W(t)\,\lambda(dt)\in[0,\infty],
        \qquad
        I_{\TT}(W)=\lim_{m\to\infty} I_{[-m,m]^l}(W).
\]
Then there exist increasing sequences \(m_k,n_k\), \(k\in\Nset\), with
\(m_k\to\infty\), \(n_k\to\infty\), and deterministic points
\[
        t_{kj}\in[-m_k,m_k]^l,\qquad 1\le j\le n_k,
\]
such that, in the extended sense in \([0,\infty]\),
\bqn{
\frac{\lambda([-m_k,m_k]^l)}{n_k}
        \sum_{j=1}^{n_k} W(t_{kj})
        \toas
        I_{\TT}(W),
        \qquad k\to\infty .
\label{Pure}
}
Moreover, the law of \(I_{\TT}(W)\) depends only on the finite-dimensional
distributions of \(W\), and
\[
        I_{\TT}(B^hW)=I_{\TT}(W),
        \qquad \mathbb P_1\text{-a.s.},\ \forall h\in\TT .
\]
\label{lemGen2}
\EEL

\def\equaldis{\stackrel{d}{=}}

\prooflem{lemGen2}
We borrow the proof idea from   \cite[Lem 10.4.2]{MR3561100}. Let
\[
        Q_m=[-m,m]^l,\qquad m\in\Nset .
\]
By joint measurability and Tonelli's theorem
\[
        I_{Q_m}(W)=\int_{Q_m}W(t)\lambda(dt)
\]
is a well-defined \([0,\infty]\)-valued rv. Moreover, we have
\[
        I_{[-m,m]^l}(W)\uparrow
        I_{\TT}(W):=\lim_{m\to\infty}I_{[-m,m]^l}(W)
        \in[0,\infty],
        \qquad \mathbb P_1\text{-a.s.}
\]

Hence, in view of \cite[Lem ~9.2]{hashorva2021shiftinvariant}, it suffices
to show that for each fixed \(m\in\Nset\) there exist points
\(t_{mj}\in[-m,m]^l\), \(j\ge1\), such that
\[
        \frac{\lambda(Q_m)}{n}\sum_{j=1}^n W(t_{mj})
        \toas
        I_{Q_m}(W),
        \qquad n\to\infty ,
\]
where the convergence is understood in \([0,\infty]\).

Fix \(m\in\Nset\), and let \(T_{mj}\), \(j\ge1\), be iid random vectors
uniformly distributed on \(Q_m\), defined on an auxiliary probability space
\((\Omega_2,\mathscr F_2,\mathbb P_2)\). We may take all arrays
\(\{T_{mj},j\ge1\}\), \(m\in\Nset\), on the same auxiliary space.

For \(\mathbb P_1\)-a.e. \(\omega_1\), the function
\[
        t\mapsto W(\omega_1,t)
\]
is Lebesgue measurable and non-negative on \(Q_m\). For such \(\omega_1\),
the rvs  
\[
        W(\omega_1,T_{mj}),\qquad j\ge1,
\]
are iid under \(\mathbb P_2\), with extended mean
\[      
        \frac{1}{\lambda(Q_m)}I_{Q_m}(W)(\omega_1)\in[0,\infty].
\]
If this extended mean is finite, the usual strong law of large numbers gives
\[
        \frac{\lambda(Q_m)}{n}
        \sum_{j=1}^n W(\omega_1,T_{mj})
        \to
        I_{Q_m}(W)(\omega_1),
        \qquad n\to\infty .
\]
If the extended mean is infinite, apply the usual strong law to
\(W(\omega_1,T_{mj})\wedge K\). Then
\[
\liminf_{n\to\infty}
        \frac{\lambda(Q_m)}{n}
        \sum_{j=1}^n W(\omega_1,T_{mj})
\ge
        \lambda(Q_m)\,
        \mathbb E_2\{W(\omega_1,T_{m1})\wedge K\}.
\]
Letting \(K\to\infty\) gives divergence to infinity, since
\[
        \lambda(Q_m)\,
        \mathbb E_2\{W(\omega_1,T_{m1})\wedge K\}
        =
        \int_{Q_m}\bigl(W(\omega_1,t)\wedge K\bigr)\lambda(dt)
        \uparrow
        I_{Q_m}(W)(\omega_1)=\infty .
\]
Thus, for \(\mathbb P_1\)-a.e. \(\omega_1\),
\[
        \frac{\lambda(Q_m)}{n}
        \sum_{j=1}^n W(\omega_1,T_{mj})
        \to
        I_{Q_m}(W)(\omega_1)
\]
in \([0,\infty]\), for \(\mathbb P_2\)-a.e. auxiliary outcome.

Hence, by Fubini, the event
\[
A_m=
\left\{
        \frac{\lambda(Q_m)}{n}
        \sum_{j=1}^n W(T_{mj})
        \to
        I_{Q_m}(W)
        \text{ in }[0,\infty]
\right\}
\]
belongs to \(\mathscr F_1\otimes\mathscr F_2\) and satisfies
\[
        (\mathbb P_1\otimes\mathbb P_2)(A_m)=1.
\]
By Fubini's theorem, there exists a set
\(\Omega_{2,m}\in\mathscr F_2\) with \(\mathbb P_2(\Omega_{2,m})=1\)
such that for every \(\omega_2\in\Omega_{2,m}\),
\[
        \frac{\lambda(Q_m)}{n}
        \sum_{j=1}^n W(T_{mj}(\omega_2))
        \toas
        I_{Q_m}(W),
        \qquad n\to\infty .
\]
Since the set of \(m\)'s is countable, choose
\[
        \omega_2^\star\in\bigcap_{m=1}^\infty \Omega_{2,m}.
\]
Setting 
\[
        t_{mj}=T_{mj}(\omega_2^\star),
        \qquad j\ge1,\ m\in\Nset 
\]
for every fixed \(m\) we obtain 
\[
        \frac{\lambda(Q_m)}{n}\sum_{j=1}^n W(t_{mj})
        \toas
        I_{Q_m}(W),
        \qquad n\to\infty 
\]
and hence  \eqref{Pure} follows from \cite[Lem 9.2]{hashorva2021shiftinvariant}.

Repeating the preceding random sampling argument simultaneously for \(W\)
and another non-negative jointly measurable RF \(W'\) with the same
finite-dimensional distributions as \(W\), and intersecting the corresponding
full \(\mathbb P_2\)-probability sets, we may choose the same deterministic
points in the finite sums for both RFs. For every fixed \(m\) and \(n\) we have 
\[
        \frac{\lambda(Q_m)}{n}\sum_{j=1}^n W(t_{mj})
        \equaldis
        \frac{\lambda(Q_m)}{n}\sum_{j=1}^n W'(t_{mj}),
\]
because \(W\) and \(W'\) have the same finite-dimensional distributions.
Passing to the almost sure limits gives
\[
        I_{Q_m}(W)\equaldis I_{Q_m}(W'),
        \qquad m\in\Nset .
\]
Letting \(m\to\infty\), and using monotone convergence, yields
\[
        I_{\TT}(W)\equaldis I_{\TT}(W').
\]
Hence the law of \(I_{\TT}(W)\) depends only on the finite-dimensional
distributions of \(W\).

Finally, for \(h\in\TT\), using the translation-invariance of Lebesgue
measure
\[
\begin{aligned}
        I_{\TT}(B^hW)
        &=
        \lim_{m\to\infty}\int_{Q_m} W(t-h)\lambda(dt)  =
        \lim_{m\to\infty}\int_{Q_m-h} W(s)\lambda(ds).
\end{aligned}
\]
Since \(Q_m-h\uparrow\TT\), monotone convergence gives
\[
        \lim_{m\to\infty}\int_{Q_m-h} W(s)\lambda(ds)
        =
        \int_{\TT}W(s)\lambda(ds)
        =
        I_{\TT}(W).
\]
This identity holds pathwise for every \(h\in\TT\) implying  
\[
        I_{\TT}(B^hW)=I_{\TT}(W),
        \qquad \mathbb P_1\text{-a.s.},\ \forall h\in\TT .
\]
\QED

\bibliographystyle{ieeetr}
\bibliography{EEEA}

\newcommand{\nosort}[1]{}\def\polhk#1{\setbox0=\hbox{#1}{\ooalign{\hidewidth \lower1.5ex\hbox{`}\hidewidth\crcr\unhbox0}}}
\begin{thebibliography}{10}

\bibitem{hashorva2021shiftinvariant}
E.~Hashorva, ``Shift-invariant homogeneous classes of random fields,'' {\em Journal of Mathematical Analysis and Applications}, p.~128517, 2024.

\bibitem{kab2009}
Z.~Kabluchko, M.~Schlather, and L.~de~Haan, ``Stationary max-stable fields associated to negative definite functions,'' {\em Ann. Probab.}, vol.~37, pp.~2042--2065, 2009.

\bibitem{KumeE}
E.~Hashorva and A.~Kume, ``Multivariate max-stable processes and homogeneous functionals,'' {\em Statist. Probab. Lett.}, vol.~173, p.~109066, 2021.

\bibitem{deHaan}
L.~de~Haan, ``A spectral representation for max-stable processes,'' {\em Ann. Probab.}, vol.~12, no.~4, pp.~1194--1204, 1984.

\bibitem{dom2016}
C.~Dombry and Z.~Kabluchko, ``Ergodic decompositions of stationary max-stable processes in terms of their spectral functions,'' {\em Stochastic Processes and their Applications}, vol.~127, no.~6, pp.~1763--1784, 2017.

\bibitem{MolchanovBE}
I.~Molchanov, M.~Schmutz, and K.~Stucki, ``Invariance properties of random vectors and stochastic processes based on the zonoid concept,'' {\em Bernoulli}, vol.~20, no.~3, pp.~1210--1233, 2014.

\bibitem{Htilt}
E.~Hashorva, ``Representations of max-stable processes via exponential tilting,'' {\em Stochastic Process. Appl.}, vol.~128, no.~9, pp.~2952--2978, 2018.

\bibitem{MartinE}
M.~Bladt, E.~Hashorva, and G.~Shevchenko, ``Tail measures and regular variation,'' {\em Electron. J. Probab.}, vol.~27, pp.~Paper No. 64, 43, 2022.

\bibitem{kulik:soulier:2020}
R.~Kulik and P.~Soulier, {\em {Heavy tailed time series.}}
\newblock Cham: Springer, 2020.

\bibitem{PH2020}
P.~Soulier, ``The tail process and tail measure of continuous time regularly varying stochastic processes,'' {\em Extremes}, vol.~25, no.~1, pp.~107--173, 2022.

\bibitem{hashorva2025cluster}
E.~Hashorva, ``Cluster random fields and random-shift representations,'' {\em Journal of Theoretical Probability}, vol.~38, no.~3, p.~50, 2025.

\bibitem{Hrovje}
H.~Planini\'{c} and P.~Soulier, ``The tail process revisited,'' {\em Extremes}, vol.~21, no.~4, pp.~551--579, 2018.

\bibitem{klem}
C.~Dombry, E.~Hashorva, and P.~Soulier, ``Tail measure and spectral tail process of regularly varying time series,'' {\em Ann. Appl. Probab.}, vol.~28, no.~6, pp.~3884--3921, 2018.

\bibitem{BojanS}
B.~Basrak and J.~Segers, ``Regularly varying multivariate time series,'' {\em Stochastic Process. Appl.}, vol.~119, no.~4, pp.~1055--1080, 2009.

\bibitem{BojanPhilippe}
B.~Basrak, H.~Planinic, and P.~Soulier, ``An invariance principle for sums and record times of regularly varying stationary sequences,'' {\em Probab. Theory Relat. Fields}, vol.~172, p.~869–914, 2018.

\bibitem{Resnickart}
S.~Resnick, {\em The Art of Finding Hidden Risks: Hidden Regular Variation in the 21st Century}.
\newblock Springer Nature, 2024.

\bibitem{wao}
T.~Owada and G.~Samorodnitsky, ``Tail measures of stochastic processes of random fields with regularly varying tails. {T}echnical report,'' 2012.

\bibitem{MR3561100}
G.~Samorodnitsky, {\em Stochastic processes and long range dependence}.
\newblock Springer Series in Operations Research and Financial Engineering, Springer, Cham, 2016.

\bibitem{MR1280932}
G.~Samorodnitsky and M.~S. Taqqu, {\em Stable non-{G}aussian random processes}.
\newblock Stochastic Modeling, Chapman \& Hall, New York, 1994.
\newblock Stochastic models with infinite variance.

\bibitem{SegersEx17}
J.~Segers, Y.~Zhao, and T.~Meinguet, ``Polar decomposition of regularly varying time series in star-shaped metric spaces,'' {\em Extremes}, vol.~20, no.~3, pp.~539--566, 2017.

\bibitem{Ilya25}
B.~Basrak, N.~Milinčević, and I.~Molchanov, ``Foundations of regular variation on topological spaces,'' {\em arXiv preprint arXiv:2503.00921}, 2025.

\bibitem{Guenter}
G.~Last, ``Tail processes and tail measures: An approach via palm calculus,'' {\em Extremes}, vol.~26, no.~4, pp.~715--746, 2023.

\bibitem{DombKabA}
C.~Dombry and Z.~Kabluchko, ``Ergodic decompositions of stationary max-stable processes in terms of their spectral functions,'' {\em Stochastic Process. Appl.}, vol.~127, no.~6, pp.~1763--1784, 2017.

\bibitem{MR4280158}
Y.~Cissokho and R.~Kulik, ``Estimation of cluster functionals for regularly varying time series: sliding blocks estimators,'' {\em Electron. J. Stat.}, vol.~15, no.~1, pp.~2777--2831, 2021.

\bibitem{chen2025asymptotic}
Z.~Chen and R.~Kulik, ``Asymptotic expansions for blocks estimators: Pot framework,'' {\em Stochastic Processes and their Applications}, p.~104744, 2025.

\bibitem{chen2023limit}
Z.~Chen and R.~Kulik, ``Limit theorems for unbounded cluster functionals of regularly varying time series,'' {\em arXiv preprint arXiv:2308.16270}, 2023.

\bibitem{StoevSPA}
S.~A. Stoev, ``On the ergodicity and mixing of max-stable processes,'' {\em Stochastic Process. Appl.}, vol.~118, no.~9, pp.~1679--1705, 2008.

\bibitem{GeoSS}
K.~Bisewski, E.~Hashorva, and G.~Shevchenko, ``The harmonic mean formula for random processes,'' {\em Stochastic Analysis and Applications}, vol.~41, no.~3, pp.~591--603, 2023.

\end{thebibliography}
\end{document}